\documentclass{amsart}
\usepackage[mathcal]{eucal}
\usepackage{amssymb, graphics}
\usepackage[dvips]{graphicx}
    \newtheorem{thm}{Theorem}
    \newtheorem{cor}[thm]{Corollary}
    \newtheorem{lem}[thm]{Lemma}
    \newtheorem{prop}[thm]{Proposition}
    \theoremstyle{definition}
    
    \theoremstyle{remark}
    \newtheorem{rem}[thm]{Remark}
    \numberwithin{equation}{section}

\begin{document}

\title[The quantum Teichm\"uller Space]
    {The quantum Teichm\"uller space \\ as a
noncommutative algebraic object}

\author{Xiaobo Liu }

\address{Department of Mathematics, MC4444, 
Columbia University, 2990 Broadway, New York, NY 10027}

\email{xiaoboli@math.columbia.edu}

\thanks{ This work was partially supported by NSF
grant DMS-0103511}

\date{\today}



\begin{abstract}
We consider the quantum Teichm\"uller space of the punctured
surface introduced by Chekhov-Fock-Kashaev, and formalize it as a
noncommutative deformation of the space of algebraic functions on
the Teichm\"uller space of the surface. In order to apply it in
3-dimensional topology, we put more attention to the details
involving small surfaces.

\end{abstract}

\maketitle

Let $S$ be an oriented surface of finite topological type, with at
least one puncture. A quantization of the
    Teichm\"uller space $\mathcal T(S)$ of $S$ was
developed by L. Chekhov and V. Fock \cite{Foc,FocChe,CheFoc} and,
independently, by R.~Kashaev \cite{Kas} (see also \cite{Tes}) as
an approach to quantum gravity in $2+1$ dimensions. This is a
deformation of the $\mathrm C^*$--algebra of functions on the
usual Teichm\"uller space $\mathcal T(S)$ of $S$, depending on a
parameter $\hbar$, in such a way that the linearization of this
deformation at $\hbar=0$ corresponds to the Weil-Petersson Poisson
structure on $\mathcal T(S)$.

In this paper, we develop a slightly different version of this
quantization, which has a more algebraic flavor. It essentially is
the image under the exponential map of the quantization of
Chekhov-Fock-Kashaev. The original quantization was expressed in
terms of self-adjoint operators on Hilbert spaces and made strong
use of the holomorphic function
$$
\phi^\hbar (z) = - \frac{\pi\hbar}2 \int_{-\infty}^{\infty}
\frac{\mathrm e ^{- \mathrm i tz}}{\sinh \pi t \, \sinh \pi \hbar
t}\,dt
$$
called the quantum dilogarithm function. Our exponential version
enables us to leave the realm of analysis and to focus on the
algebraic aspects of the construction.

 From a mathematical point of view, the main benefit of the quantum
Teichm\"uller space $\mathcal T_S^q$ which we construct here is
that it admits a rich finite-dimensional representation theory.
This representation theory is investigated in \cite{BonLiu}, where
we show that it is strongly connected to 3--dimensional hyperbolic
geometry through the space of representations $\pi_1(S)
\rightarrow \mathrm{PSL}_2(\mathbb C)$. In particular, we use this
connection to construct in \cite{BonLiu} quantum invariants for
diffeomorphisms of $S$. The present paper is devoted to laying
down the foundations of the theory, by showing that the quantum
Teichm\"uller space $\mathcal T_S^q$ is well-defined as an object
in non-commutative algebraic geometry.

Our construction parallels that of Chekhov and Fock. We start from
Thurston's exponential shear coordinates for Teichm\"uller space,
associated to an ideal triangulation of $S$. A fundamental
property of these exponential shear coordinates is that, as one
shifts from one ideal triangulation to another, the corresponding
coordinate changes are rational. In particular, there is a
well-defined notion of rational function on $\mathcal T(S)$. In
general, the quantization of a space is a deformation of the
algebra of continuous functions on this space. Here we construct a
deformation of the algebra $\mathrm{Rat}\,\mathcal T(S)$ of all
rational functions on $\mathcal T(S)$. The main strategy is, first
to define a deformation of the expression of
$\mathrm{Rat}\,\mathcal T(S)$ in the set of shear coordinates
associated to an ideal triangulation, and then to construct
appropriate ``coordinate change isomorphisms'' to make the
construction independent of any choice of ideal triangulation. The
technical challenge is to check that these coordinate change
isomorphisms are compatible with each other, which is accomplished
in \S\ref{sect:CheFock}-\S\ref{sect:QuantumTeich}.

Many aspects of this construction are already implicitly or
explicitly present in the work of Chekhov-Fock. In particular, the
algebraic coordinate change isomorphisms originate from the
symmetries of the quantum dilogarithm function $\phi^\hbar$. Our
main contribution is to systematically develop the theory from an
algebraic point of view, and to carefully check the details of the
construction. From a technical point of view, this leads us to be
a little more careful with non-embedded diagonal exchanges, and to
explicitly describe the formulas which occur in these cases.

So far, we have avoided defining the Teichm\"uller space of the
surface $S$. The reason is that there are several competing
definitions, leading to different objects. One possibility is to
define the Teichm\"uller space $\mathcal T(S)$ as the space of
isotopy classes of complete hyperbolic metrics on $S$, possibly
with infinite area. This is a manifold-with-corners, of dimension
$6g-6+3p$ where $g$ is the genus of $S$ and $p$ is the number of
its punctures. This space can be expanded to the \emph{enhanced
Teichm\"uller space} $\widetilde{\mathcal T}(S)$ defined in
\S\ref{sect:Shear}, which has the same space of rational functions
as $\mathcal T(S)$. This is the set-up of \cite{FocChe} and of the
first part of this article. However, a more commonly used object
is the \emph{cusped Teichm\"uller space} $\mathcal{CT}(S) \subset
\widetilde{\mathcal T}(S)$, consisting of all isotopy classes of
finite area complete hyperbolic metrics on $S$; this is a manifold
of dimension $6g-6+2p$. We construct a quantization of this cusped
Teichm\"uller space in the last section \S\ref{sect:QuantumCusped}
of the article.

There is one property which we would like to mention before
closing this introduction. The original definition of the quantum
Teichm\"uller space was grounded in the geometry of hyperbolic
metrics on $S$. However, a recent result of Hua Bai \cite{Bai}
shows that it is intrinsically tied to the combinatorics of the
Harer-Penner simplicial complex of ideal triangulations \cite{Har,
Pen}. Indeed, Bai proved that the coordinate change isomorphisms
considered in \S\ref{sect:CheFock} are the only ones which satisfy
a certain number of natural conditions. This emphasizes the
algebraic and combinatorial nature of the quantum Teichm\"uller
space.

\medskip
\noindent{Acknowledgements:} We would like to thank Francis
Bonahon, Leonid Chekhov and Bob Penner for their help and
encouragements.

\section{Ideal Triangulations}
\label{sect:IdealTri}

Let $S$ be an oriented surface of genus $g$ with $p\geq 1$
punctures, obtained by removing $p$ points $\{v_1,\ldots,v_p\}$
from the closed oriented surface $\bar S$ of genus $g$. An
\emph{ideal triangulation} of $S$ is a triangulation of the closed
surface $\bar{S}$ whose vertex set is exactly
$\{v_1,\ldots,v_p\}$. An Euler characteristic argument shows that
any ideal triangulation has $n=6g-6+3p$ edges.

\begin{figure}[h]
\includegraphics[scale=.75]{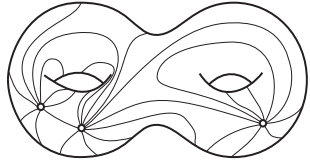}
\caption{An ideal triangulation of a genus 2  surface with 3
punctures}
\end{figure}

Two ideal triangulations are considered the same if they are
isotopic. In addition, we require that each ideal triangulation
$\lambda$ is endowed with an indexing $\lambda_1$,
$\lambda_2$,\dots, $\lambda_n$ of its edges by the index set
$\{1,2, \dots, n\}$. Let $\Lambda(S)$ denote the set of isotopy
classes of such indexed ideal triangulations $\lambda$.

The set $\Lambda(S)$ admits a natural action of the group
$\mathfrak{S}_n$ of permutations of $n$ elements, acting by
permuting the indices of the edges of $\lambda$. Namely $\lambda'
= \alpha(\lambda)$ for $\alpha \in\mathfrak{S}_n$ if its $i$--th
edge $\lambda'_i$ is equal to $\lambda_{\alpha(i)}$.

\begin{figure}[h]
\includegraphics{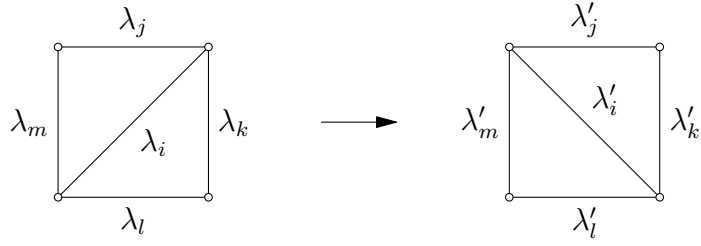}
\caption{A diagonal exchange} \label{fig:DiagExch}
\end{figure}

Another important transformation of $\Lambda(S)$ is provided by
the \emph {$i$--th diagonal exchange map} $\Delta_i:
\Lambda(S)\rightarrow \Lambda(S)$ defined as follows. The $i$--th
edge $\lambda_i$ of an ideal triangulation $\lambda\in \Lambda(S)$
is adjacent to two triangles. If these two triangles are distinct,
their union forms a square $Q$ with diagonal $\lambda_i$. Then
$\Delta_i(\lambda)$ is obtained from $\lambda$ by replacing the
edge $\lambda_i$ by the other diagonal $\lambda_i'$ of the square
$Q$, as illustrated in Figure~\ref{fig:DiagExch}. By convention,
$\Delta_i(\lambda)=\lambda$ when the two sides of $\lambda_i$
belong to the same triangle; this happens exactly when $\lambda_i$
is the only edge of $\lambda$ leading to a certain  puncture of
$S$.

\begin{figure}[h]

\includegraphics[scale=0.5]{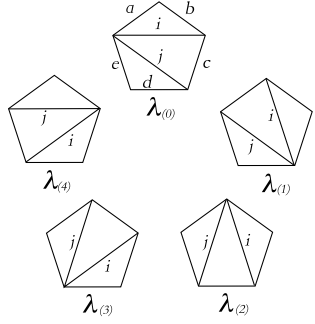}
\caption{The Pentagon Relation} \label{fig:pentagon}
\end{figure}

The  reindexings and diagonal exchanges satisfy the following
relations:
\begin{enumerate}

\item The \emph{Composition Relation}: $(\alpha\beta)(\lambda) =
\alpha(\beta(\lambda))$ for every $\alpha$,
$\beta\in\mathfrak{S}_n$;

\item The \emph{Reflexivity Relation}: $(\Delta_i)^2=\mathrm{Id}$;

\item The \emph{Reindexing Relation}: $\Delta_i \circ \alpha  =
\alpha \circ \Delta_{\alpha(i)} $ for every
$\alpha\in\mathfrak{S}_n$;

\item The \emph{Distant Commutativity Relation}: If $\lambda_i$
and $\lambda_j$ do not belong to a same triangle of $\lambda\in
\Lambda(S)$, then $\Delta_i \circ \Delta_j(\lambda)= \Delta_j
\circ \Delta_i(\lambda)$;

\item The \emph{Pentagon Relation}: If three triangles of an ideal
triangulation $\lambda\in \Lambda(S)$ form a pentagon with
diagonals $\lambda_i$, $\lambda_j$, and if
$\alpha_{i\leftrightarrow j} \in\mathfrak{S}_n$ denotes the
transposition exchanging $i$ and $j$, then
\begin{equation}
\label{eqn:pentagon} \Delta_i\circ \Delta_j\circ \Delta_i\circ
\Delta_j\circ \Delta_i(\lambda) = \alpha_{i\leftrightarrow
j}(\lambda).
\end{equation}

The Pentagon Relation is illustrated in Figure~\ref{fig:pentagon}.

\end{enumerate}

The following two results of R.~Penner \cite{Pen} (although the
first one may have an earlier history) are crucial in the
construction of the quantum Teichm\"uller space. We refer to
\cite{Pen} for their proof.

\begin{thm}
\label{thm:Penner1} Given  two ideal triangulations
$\lambda,\lambda'\in \Lambda(S)$, there exists a finite sequence
of ideal triangulations $\lambda=\lambda_{(0)}$, $\lambda_{(1)}$,
\dots, $\lambda_{(m)}=\lambda'$ such that each $\lambda_{(k+1)}$
is obtained from $\lambda_{(k)}$ by a diagonal exchange or by a
reindexing of its edges. \qed
\end{thm}

\begin{thm}
\label{thm:Penner2} Given  two ideal triangulations
$\lambda,\lambda'\in \Lambda(S)$ and given two sequences
$\lambda=\lambda_{(0)}$, $\lambda_{(1)}$, \ldots, $\lambda_{(m)}
=\lambda'$ and $\lambda=\lambda'_{(0)}$, $\lambda'_{(1)}$, \ldots,
$\lambda'_{(m')} =\lambda'$ of diagonal exchanges and reindexings
connecting them as in Theorem~\ref{thm:Penner1}, these two
sequences can be related to each other by successive applications
of the following moves and of their inverses.
\begin{enumerate}

\item Use the Composition Relation to replace \dots,
$\lambda_{(k)}$, $\beta(\lambda_{(k)})$, $\alpha\circ
\beta(\lambda_{(k)})$, \dots\ by \dots, $\lambda_{(k)}$, $(\alpha
\beta)(\lambda_{(k)})$, \dots\ where $\alpha$, $\beta\in
\mathfrak{S}_n$.

\item Use the Reflexivity Relation to replace \dots,
$\lambda_{(k)}$, \dots\ by \dots, $ \lambda_{(k)}$,
$\Delta_i(\lambda_{(k)}) $, $\lambda_{(k)}$, \dots\ .

\item Use the Reindexing Relation to replace \dots,
$\lambda_{(k)}$, $\alpha(\lambda_{(k)})$, $\Delta_i \circ
\alpha(\lambda_{(k)}) =\lambda_{(k+2)}$, \dots\ by \dots ,
$\lambda_{(k)}$, $\Delta_{\alpha(i)}(\lambda_{(k)})$, $\alpha\circ
\Delta_{\alpha(i)} (\lambda_{(k)})=\lambda_{(k+2)}$, \dots\ where
$\alpha\in \mathfrak{S}_n$.

\item Use the Distant Commutativity Relation to replace  \dots,
$\lambda_{(k)}$, \dots\ by \dots, $ \lambda_{(k)}$,
$\Delta_i(\lambda_{(k)})$, $\Delta_j\Delta_i(\lambda_{(k)})$,
$\Delta_j(\lambda_{(k)})$, $\lambda_{(k)}$, \dots\  where
$\lambda_i,\lambda_j$ are two edges which do not belong to a same
triangle of $\lambda_{(k)}$.

\item Use the Pentagon Relation to replace \dots, $\lambda_{(k)}$,
\dots\ by \dots $\lambda_{(k)}$, $\Delta_i(\lambda_{(k)})$,
$\Delta_j\circ \Delta_i(\lambda_{(k)})$, $\Delta_i\circ
\Delta_j\circ \Delta_i(\lambda_{(k)})$, $\Delta_j\circ
\Delta_i\circ \Delta_j\circ \Delta_i(\lambda_{(k)})$,
$\alpha_{i\leftrightarrow j} (\lambda_{(k)})$, $\lambda_{(k)}$,
\dots\ where $\lambda_i,\lambda_j$ are two diagonals of a pentagon
of $\lambda_{(k)}$.
    \qed

\end{enumerate}
\end{thm}

\section{The exponential shear coordinates for
Teichm\"uller space} \label{sect:Shear}

W.~Thurston associated to an ideal triangulation $\lambda$ a
certain system of \emph{shear coordinates} for the Teichm\"uller
space $\mathcal T(S)$. See \cite{Thu86} for the dual notion of
length coordinates, and \cite{Bon}\cite{Foc}\cite{FocChe} for
details on these shear coordinates.

Consider a complete hyperbolic metric $m \in \mathcal{T}(S)$. It
is well-known that the ends of the complete hyperbolic surface
$(S, m)$ can be of two types: finite area \emph{cusps} bounded on
one side by a horocycle; and  infinite area \emph{funnels} bounded
on one side by a simple closed geodesic. This can also be
expressed in terms of the \emph{convex core} $\mathrm{Conv}(S,m)$
of $(S,m)$, which is the smallest non-empty closed convex subset
of $(S,m)$, and is bounded in $S$ by a family of disjoint simple
closed geodesics. The cusp ends of $(S,m)$ are those which are
also ends of $\mathrm{Conv}(S,m)$, while each funnel end of $S$
faces a boundary component of $\mathrm{Conv}(S,m)$. Note that the
interior $\mathrm{Int} (\mathrm{Conv}(S,m))$ is homeomorphic to
$S$, by the homeomorphism which is uniquely determined up to
isotopy by the property that it is homotopic to the inclusion map.

It is convenient to enhance the hyperbolic metric $m\in \mathcal
T(S)$ with some additional data, consisting of an orientation for
the boundary $\partial\mathrm {Conv}(S,m)$ of its convex cores.
Let the \emph{enhanced Teichm\"uller space} $\widetilde{\mathcal
T}(S)$ consist of all isotopy classes of hyperbolic metrics $m\in
\mathcal T(S)$ enhanced with an orientation of $\partial\mathrm
{Conv}(S,m)$. Since the convex core $\mathrm {Conv}(S,m)$ depends
continuously on the metric $m$, the enhanced Teichm\"uller space
$\widetilde{\mathcal T}(S)$ inherits from the topology of
$\mathcal T(S)$ a topology for which the natural projection
$\widetilde{\mathcal T}(S) \rightarrow \mathcal T(S)$ is a
branched covering map.

\begin{figure}[h]
\includegraphics[scale=0.75]{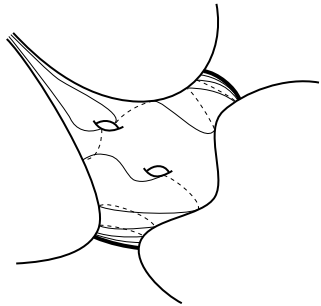}
\caption{A hyperbolic surface} \label{fig:hypsurf}
\end{figure}

Consider an enhanced hyperbolic metric $m \in \widetilde{\mathcal
T}(S)$ together with an ideal triangulation $\lambda$. Each edge
$\lambda_i$ specifies a proper homotopy class of paths going from
one end of $\mathrm{Int} (\mathrm{Conv}(S,m)) \cong S$ to another
end. This proper homotopy class is also realized by a unique
$m$--geodesic $g_i$ such that each end of $g_i$, either converges
towards a cusp end of $S$,  or spirals around a component of
$\partial \mathrm{Conv}(S,m)$ in the direction specified by the
enhancement of $m$. The union of $\partial \mathrm{Conv}(S,m)$ and
of the $g_i$ forms an $m$--geodesic lamination $g$ contained in
the convex core $\mathrm{Conv}(S,m)$.

The enhanced hyperbolic metric $m \in \widetilde{\mathcal T}(S)$
now associates to the edge $\lambda_i$ of $\lambda$ a positive
number $x_i$ defined as follows. The geodesic $g_i$ separates two
triangle components $T_i^1$ and $T_i^2$ of $\mathrm{Conv}(S,m)-g$.
Identify the universal covering of S, endowed with the metric $m$,
to the hyperbolic plane $\mathbb H^2$. Lift $g_i$, $T_i^1$ and
$T_i^2$ to a geodesic $\widetilde g_i$ and two triangles
$\widetilde T_i^1$ and $\widetilde T_i^2$ in  $\mathbb H^2$ so
that the union $\widetilde g_i \cup \widetilde T_i^1\cup
\widetilde T_i^2$ forms a square $\widetilde Q$ in $\mathbb H^2$.
In the standard upper half-space model for $\mathbb H^2$, let
$z_-$, $z_+$, $z_{\mathrm r}$, $z_{\mathrm l}$ be the vertices of
$\widetilde Q$ in such a way that $\widetilde g_i$ goes from $z_-$
to $z_+$ and, for this orientation of $\widetilde g_i$,
$z_{\mathrm r}$, $z_{\mathrm l}$ are respectively to the right and
to the left of $\widetilde g_i$ for the orientation of $\widetilde
Q$ given by the orientation of $S$. Then,
$$x_i=- \,\textrm{cross-ratio}\,
(z_{\mathrm r},z_{\mathrm l},z_-,z_+)=-\, \frac{(z_{\mathrm
r}-z_-)(z_{\mathrm l}-z_+)}{(z_{\mathrm r}-z_+) (z_{\mathrm
l}-z_-)}.$$ Note that $x_i$ is positive since the points $z_-$,
$z_{\mathrm r}$, $z_+$, $z_{\mathrm l}$ occur in this order in the
real line bounding the upper half-space $\mathbb H^2$.

The real numbers $x_i$ are the \emph{exponential shear
coordinates} of the enhanced hyperbolic metric $m \in
\widetilde{\mathcal T}(S)$. The more standard shear coordinates
are just their logarithms $\log x_i$, but the $x_i$ turn out to be
better behaved for our purposes.

There is an inverse construction which associates a hyperbolic
metric to each system of positive weights $x_i$ attached to the
edges $\lambda_i$ of the ideal triangulation $\lambda$: Identify
each of the components of $S-\lambda$ to a triangle with vertices
at infinity in $\mathbb H^2$, and glue these hyperbolic triangles
together in such a way that adjacent triangles form a square whose
vertices have cross-ration $-x_i$ as above. This defines a
\emph{possibly incomplete} hyperbolic metric on the surface $S$.
An analysis of this metric near the ends of $S$ shows that its
completion is a hyperbolic surface $S'$ with geodesic boundary,
and that each end of an edge of $\lambda$ either spirals towards a
component of $\partial S'$ or converges towards a cusp end of
$S'$. Extending $S'$ to a complete hyperbolic surface without
boundary and identifying this surface to $S$, we now have found a
complete hyperbolic metric $m$ on $S$ whose convex core is
isometric to $S'$. In addition, the spiralling pattern of the ends
of $\lambda$ provides an enhancement for the hyperbolic metric
$m$. See \cite[\S 3.4]{Thu97} or \cite{Bon} for details.

The $x_i$ then define a homeomorphism $\phi_\lambda
:\widetilde{\mathcal T}(S) \rightarrow \mathbb R_+^n$ between the
enhanced Teichm\"uller space $\widetilde{\mathcal T}(S)$ and
$\mathbb R_+^n$.

For the enhanced hyperbolic metric $m \in \widetilde{\mathcal
T}(S)$ associated to the exponential shear parameters $x_i$, the
geometry of $m$ near the $j$--th puncture $v_j$ is completely
determined by the \emph{exponential length parameter}
$p_j=x_1^{k_{1j}}x_2^{k_{2j}}\dots x_n^{k_{nj}}$ where $k_{ij}\in
\{0,1,2\}$ is the number of end points of the edge $\lambda_i$
that are equal to $v_j$. The end of $S$ corresdonding to $v_j$ is
a cusp for $m$ exactly when $p_j=1$. Otherwise, this end faces a
boundary component of $\mathrm{Conv}(S,m)$ which is a closed
geodesic of length $\lvert \log p_j \rvert$, and the orientation
of this boundary component defined by the enhancement of $m$
coincides with the boundary orientation exactly when $p_j<1$.

The exponential shear  coordinates associates a parametrization
$\phi_\lambda :\widetilde{\mathcal T}(S) \rightarrow \mathbb
R_+^n$ to each ideal triangulation $\lambda\in \Lambda(S)$
(endowed with an indexing of its edges). We now investigate the
coordinate changes $\phi_{\lambda'} \circ \phi_{\lambda}^{-1}$
associated to two ideal triangulations.

If $\lambda'= \alpha(\lambda)$ is obtained by reindexing the edges
of $\lambda$ by $\alpha \in \mathfrak S_n$, it is immediate that
$\phi_{\lambda'} \circ \phi_{\lambda}^{-1}$ is just the
permutation of the coordinates by $\alpha$. For a diagonal
exchange, the expression for $\phi_{\lambda'} \circ
\phi_{\lambda}^{-1}$ depends on the possible identifications
between the sides of the square where the diagonal exchange takes
place.

\begin{prop}
\label{prop:DiagExch} Suppose that the ideal triangulations
$\lambda$, $\lambda'\in \Lambda(S)$ are obtained from each other
by a diagonal exchange, namely that $\lambda' =
\Delta_i(\lambda)$. Label the edges of $\lambda$ involved in this
diagonal exchange as $\lambda_i$, $\lambda_j$, $\lambda_k$,
$\lambda_l$, $\lambda_m$ as in Figure~\ref {fig:DiagExch}. If
$(x_1, x_2, \dots, x_n)=\phi_\lambda(m)$ and $(x_1', x_2', \dots,
x_n')=\phi_{\lambda'}(m)$ are the exponential shear coordinates
associated to the same enhanced hyperbolic metric $m \in
\widetilde{\mathcal T}(S)$, then $x_h'=x_h$ for every $h \not\in
\{i,j,k,l,m\}$, $x_i' = x_i^{-1}$ and:
\begin{description}
\item[Case 1] if the edges $\lambda_j$, $\lambda_k$, $\lambda_l$,
$\lambda_m$  are distinct, then
$$ x'_j = (1+x_i)x_j\quad x'_k  =
(1+x_i^{-1})^{-1}x_k \quad x'_l = (1+x_i)x_l \quad x'_m =
(1+x_i^{-1})^{-1}x_m;$$

\item[Case 2] if $\lambda_j$ is identified with $\lambda_k$, and
$\lambda_l$ is distinct from $\lambda_m$, then
$$
    x'_j  =  x_ix_j \quad x'_l  =
(1+x_i)x_l \quad x'_m  =  (1+x_i^{-1})^{-1}x_m;
$$

\item[Case 3] (the inverse of Case 2) if $\lambda_j$ is identified
with $\lambda_m$, and $\lambda_k$ is distinct from $\lambda_l$,
then
$$
    x'_j  =  x_ix_j
\quad x'_k  =  (1+x_i^{-1})^{-1}x_k \quad x'_l  = (1+x_i)x_l;
$$

\item[Case 4] if $\lambda_j$ is identified with $\lambda_l$, and
$\lambda_k$ is distinct from $\lambda_m$, then
$$
    x'_j  =  (1+x_i)^2x_j
\quad x'_k  =  (1+x_i^{-1})^{-1}x_k \quad x'_m  =
(1+x_i^{-1})^{-1}x_m
$$

\item[Case 5] (the inverse of Case 4) if $\lambda_k$ is identified
with $\lambda_m$, and $\lambda_j$ is distinct from $\lambda_l$,
then
$$ x'_j  =  (1+x_i)x_j \quad x'_k  =
(1+x_i^{-1})^{-2}x_k \quad x'_l = (1+x_i)x_l;
$$

\item[Case 6] if $\lambda_j$ is identified with $\lambda_k$, and
$\lambda_l$ is identified with $\lambda_m$ (in which case $S$ is a
$3$--times punctured sphere), then
$$
    x'_j  =  x_ix_j \quad x'_l  =
x_ix_l;
$$

\item[Case 7] (the inverse of Case 6) if $\lambda_j$ is identified
with $\lambda_m$, and $\lambda_k$ is identified with $\lambda_l$
(in which case $S$ is a $3$--times punctured sphere), then
$$
    x'_j  =  x_ix_j
\quad x'_k  =  x_ix_k;
$$

\item[Case 8] if $\lambda_j$ is identified with $\lambda_l$, and
$\lambda_k$ is identified with $\lambda_m$ (in which case $S$ is a
once punctured torus), then
$$ x'_j  =  (1+x_i)^2x_j \quad x'_k  =
(1+x_i^{-1})^{-2}x_k.
$$

\end{description}
\end{prop}

\begin{proof} This immediately follows from
the combinatorics of cross-ratios. For instance, let us focus on
Case~4. The other cases are similar.

Isometrically identify $S$, endowed with the metric $m$, to the
quotient of the upper half-plane model for $\mathbb{H}^2$ under an
action of the fundamental group $\pi_1(S)$. Lift the square $Q
\subset S$ bounded by $\lambda_j$, $\lambda_k$,
$\lambda_l=\lambda_j$ and $\lambda_m$ to a square $\widetilde Q
\subset \mathbb{H}^2$ bounded by lifts $\widetilde\lambda_j$,
$\widetilde\lambda_k$, $\widetilde\lambda_l$ and
$\widetilde\lambda_m$, and with diagonal $\widetilde\lambda_i$.
The fact that $\lambda_l=\lambda_j$ means that there exists a
covering translation $\gamma \in \pi_1(S)$ such that
$\gamma\widetilde\lambda_l= \widetilde\lambda_j$. Label the
corners of the square $\widetilde Q$ clockwise as $z_-$,
$z_{\mathrm l}$, $z_+$, $z_{\mathrm r}$, starting from the corner
$\lambda_j\cap\lambda_k$. Then, by definition,
\begin{equation*}
x_i = -\frac{(z_{\mathrm r}-z_-)(z_{\mathrm l}-z_+) }{(z_{\mathrm
r}-z_+)(z_{\mathrm l}-z_-) }
\end{equation*}
   and
\begin{equation*}
x_i' = -\frac{(z_- -z_{\mathrm l})(z_+-z_{\mathrm r}) }{(z_-
-z_{\mathrm r})(z_+ - z_{\mathrm l}) } =x_i^{-1}
\end{equation*}
since the component $\widetilde\lambda_i'$ goes from $z_{\mathrm
l}$ to $z_{\mathrm r}$.

The geodesic $\widetilde\lambda_k$  is the diagonal of a square
with vertices $z_+$, $z_-$, $z_{\mathrm l}$ and a fourth vertex
$z_k$ outside of $\widetilde Q$. Then
\begin{equation*}
x_k = -\frac{(z_{k}-z_{\mathrm l}) (z_+ -z_-)}{(z_k -z_-) (z_+
-z_{\mathrm l})}.
\end{equation*}
With respect to the preimage $\widetilde\lambda' \subset \mathbb
H^2$ of $\lambda'$, $\widetilde\lambda_k' =\widetilde\lambda_k$ is
the diagonal of a square with vertices $z_-$, $z_{\mathrm r}$,
$z_{\mathrm l}$ and $z_k$. Consequently,
\begin{equation*}
x_k' = -\frac{(z_{k}-z_{\mathrm l}) (z_{\mathrm r} -z_-)}{(z_k
-z_-) (z_{\mathrm r} -z_{\mathrm l})}=(1+x_i^{-1})^{-1}x_k.
\end{equation*}
The same argument shows that $x_m'=(1+x_i^{-1})^{-1}x_m$.

For $x_j$, there is a new twist because $\gamma\in \pi_1(S)$ sends
$\widetilde \lambda_l$ to $\widetilde\lambda_j$, so that
$\widetilde\lambda_j$ is the diagonal of a square with vertices
$z_-=\gamma z_{\mathrm l}$, $z_+$, $z_{\mathrm r}=\gamma z_+$ and
$\gamma z_-$. Switching from $\lambda$ to $\lambda'$ now changes
two vertices of this square in the sense that, with respect to
$\widetilde\lambda'$, $\widetilde\lambda_j' =\widetilde\lambda_j$
is the diagonal of a square with vertices $z_-\gamma z_{\mathrm
l}$, $z_{\mathrm l}$, $z_{\mathrm r}\gamma z_+$ and $\gamma
z_{\mathrm r}$. Therefore,
\begin{equation*}
x_j = -\frac {(\gamma z_- - \gamma z_{\mathrm l}) ( z_+
-z_{\mathrm r})} {( \gamma z_- - \gamma z_+ ) ( z_+ - z_- )}
\end{equation*}
and, since $\rho(\gamma) \in \mathrm{Isom}^+(\mathbb H^3)$
respects cross-ratios,
\begin{equation*}
\begin{split}
x_j' &= -\frac {(\gamma z_{\mathrm r} -\gamma z_{\mathrm l}) (
z_{\mathrm l} -z_{\mathrm r} )} {(\gamma z_{\mathrm r} - \gamma
z_+ ) ( z_{\mathrm l} - z_-
)}\\
&=\frac {(\gamma z_{\mathrm r} -\gamma z_{\mathrm l} ) (\gamma z_-
- \gamma z_+ ) } {(\gamma z_{\mathrm r} - \gamma z_+ ) (\gamma z_-
- \gamma z_{\mathrm l} )}\, \frac { ( z_{\mathrm l} -z_{\mathrm r}
)
   ( z_+ - z_-
   )} {  ( z_{\mathrm l} - z_-  )
   ( z_+ -z_{\mathrm r}  )}\, x_j\\
&=\left(\frac { ( z_{\mathrm l} -z_{\mathrm r}  )
   ( z_+ - z_- )} {  ( z_{\mathrm l} - z_-
   )  ( z_+ -z_{\mathrm r}  )}\right)^2 x_j =
\bigl (1 +x_i \bigr)^2 x_j.
\end{split}
\end{equation*}
The fact that $x_h'=x_h$ for $h\not\in\{i,j,k,l,m\}$ is obvious
from  the definition of the shear coordinates.
\end{proof}

Using Theorem~\ref{thm:Penner1}, an immediate corollary of
Proposition~\ref{prop:DiagExch} is the following.

\begin{cor}
For two ideal triangulations $\lambda$, $\lambda' \in \Lambda(S)$,
the coordinate change map $\phi_{\lambda'} \circ
\phi_{\lambda}^{-1}: \mathbb R^n_+ \rightarrow \mathbb R^n_+$ is
rational. \qed
\end{cor}

In particular, there is a well-defined notion of rational function
    on
$\widetilde{\mathcal T}(S)$. Namely, the (partially defined)
function $f: \widetilde{\mathcal T}(S) \rightarrow \mathbb C$ is
\emph{rational} if, for an arbitrary ideal triangulation
$\lambda$, $f\circ \phi_\lambda^{-1} : \mathbb R^n_+\rightarrow
\mathbb C$ is a rational function in the usual sense. If
$\mathrm{Rat}\,\widetilde{\mathcal T}(S)$ denotes the algebra of
such rational functions, every ideal triangulation $\lambda$
specifies an algebra isomorphism $\Phi_\lambda$ from
$\mathrm{Rat}\,\widetilde{\mathcal T}(S)$ to the algebra of
rational fractions $\mathbb C(X_1, X_2, \dots, X_n)$. In
particular, at the rational function algebra level, two ideal
triangulations $\lambda$, $\lambda' \in \Lambda(S)$ induce a
coordinate change isomorphism
\begin{equation}
\label{eqn:NonQuantCoordCh} \Phi_{\lambda\lambda'}: \mathbb C(X_1,
X_2, \dots, X_n) \rightarrow \mathbb C(X_1, X_2, \dots, X_n)
\end{equation}
defined by $g\mapsto g\circ \phi_{\lambda'} \circ
\phi_{\lambda}^{-1}$.

\section{The Chekhov-Fock algebra}
\label{sect:CheFock}

The goal of this section is to quantize the enhanced Teichm\"uller
space $\widetilde{\mathcal{T}}(S)$ by defining a deformation of
the algebra $\mathrm{Rat}\,\widetilde{\mathcal T}(S)$, depending
on a parameter $q$.

Fix an ideal triangulation $\lambda\in \Lambda(S)$. The complement
$S-\lambda$ has $2n$ spikes converging towards the punctures, and
each spike is delimited by one $\lambda_i$ on one side and one
$\lambda_j$ on the other side, with possibly $i=j$. For $i$, $j\in
\{1, \dots, n\}$, let $a^{\lambda}_{ij}$ denote the number of
spikes of $S-\lambda$ which are delimited on the left by
$\lambda_i$ and on the right by $\lambda_j$, and set
\begin{equation*}
\sigma^{\lambda}_{ij}= a^{\lambda}_{ij}-a^{\lambda}_{ji}.
\end{equation*} Note that
$\sigma_{ij}^\lambda$ can only belong to the set $\{-2, -1, 0, +1,
+2\}$, and that $\sigma_{ji}^\lambda= -\sigma_{ij}^\lambda$. It
turns out that the coefficients $\sigma_{ij}^\lambda$ are related
to the expression of the Weil-Petersson symplectic form of
$\widetilde{\mathcal{T}}(S)$ in the exponential shear coordinates,
but this plays only a historical r\^ole in this paper.

The \emph{Chekhov-Fock algebra} associated to the ideal
triangulation $\lambda$ is the algebra $\mathcal{T}^q_{\lambda}$
defined by generators $X_1$, $X_1^{-1}$, $X_2$, $X_2^{-1}$, \dots,
$X_n$, $X_n^{-1}$, with each pair $X_i^{\pm 1}$ associated to an
edge $\lambda_i$ of $\lambda$, and by the relations
\begin{gather}
X_iX_j=q^{2\sigma^{\lambda}_{ij}}X_jX_i
\label{eqn:SkewCom}\\
X_iX_i^{-1}=X_i^{-1}X_i=\mathbf{1}.
\end{gather}

This is an iterated skew Laurent polynomial algebra. In
particular, it is a Noetherian ring and a right Ore domain, so
that we can introduce its \emph{fraction division algebra}
$\widehat{\mathcal{T}}^q_{\lambda}$. See for instance \cite{Coh,
GooWar, BroGoo}. The algebra $\widehat{\mathcal{T}}^q_{\lambda}$
consists of all formal fractions $PQ^{-1}$ with $P$, $Q\in
\mathcal{T}^q_{\lambda}$ and $Q\neq 0$, and two such fractions
$P_1Q_1^{-1}$ and $P_2Q_2^{-1}$ are identified if there exists
$S_1$, $S_2\in \mathcal{T}^q_{\lambda} - \{0\}$ such that $P_1S_1
= P_2S_2$ and $Q_1S_1 = Q_2S_2$.

In practice,  the Chekhov-Fock algebra ${\mathcal T}_\lambda^q$
consists of all Laurent polynomials in variables  $X_1$, $X_2$,
\dots, $X_n$ satisfying the skew-commutativity relations
(\ref{eqn:SkewCom}). Its fraction division algebra
$\widehat{\mathcal{T}}^q_{\lambda}$ consists of all rational
fractions in the variables $X_1$, $X_2$, \dots, $X_n$ satisfying
the same skew-commutativity relations. In particular, when $q=1$,
${\mathcal T}_\lambda^1$ and $\widehat{\mathcal{T}}^1_{\lambda}$
respectively coincide with the Laurent polynomial algebra $\mathbb
C [X_1^{\pm1}, X_2^{\pm1}, \dots, X_n^{\pm1}]$ and the rational
fraction algebra $\mathbb C (X_1, X_2, \dots, X_n)$. The general
${\mathcal T}_\lambda^q$ and $\widehat{\mathcal{T}}^q_{\lambda}$
can be considered as deformations of this case.

The algebras ${\mathcal T}_\lambda^q$ and
$\widehat{\mathcal{T}}^q_{\lambda}$ strongly depend on the ideal
triangulation $\lambda$. To define a triangulation independent
deformation of the algebra $\mathrm{Rat}\,\widetilde{\mathcal
T}(S)$, we need to generalize the coordinate change isomorphism of
(\ref{eqn:NonQuantCoordCh}) to this non-commutative context, by
introducing appropriate algebra isomorphisms $
\Phi_{\lambda\lambda'}^q: \widehat{\mathcal{T}}^q_{\lambda'}
\rightarrow \widehat{\mathcal{T}}^q_{\lambda} $.

There is no geometry to guide us here, so we will do this
stepwise, using Theorems~\ref{thm:Penner1} and \ref{thm:Penner2}.
To ease the exposition, we will denote by $X'_1$, $X_2'$, \dots,
$X_n'$ the generators of $\widehat{\mathcal{T}}^q_{\lambda'}$
associated to the edges $\lambda_1'$, $\lambda_2'$, \dots,
$\lambda_n'$ of $\lambda'$, and by $X_1$, $X_2$, \dots, $X_n$ the
generators of $\widehat{\mathcal{T}}^q_{\lambda}$ associated to
the edges $\lambda_1$, $\lambda_2$, \dots, $\lambda_n$ of
$\lambda$. In particular, the isomorphism
$\Phi_{\lambda\lambda'}^q: \widehat{\mathcal{T}}^q_{\lambda'}
\rightarrow \widehat{\mathcal{T}}^q_{\lambda}$ will be completely
determined once we specify the images
$\Phi_{\lambda\lambda'}^q(X_i')$.

\begin{prop}
\label{prop:QuantumDiagExch} Suppose that the ideal triangulations
$\lambda$, $\lambda'\in \Lambda(S)$ are obtained from each other
by a diagonal exchange, namely that $\lambda' =
\Delta_i(\lambda)$. Label the edges of $\lambda$ involved in this
diagonal exchange as $\lambda_i$, $\lambda_j$, $\lambda_k$,
$\lambda_l$, $\lambda_m$ as in Figure~\ref {fig:DiagExch}. Then
there is a unique algebra isomorphism
\begin{equation*}
\Phi_{\lambda\lambda'}^q: \widehat{\mathcal{T}}^q_{\lambda'}
\rightarrow \widehat{\mathcal{T}}^q_{\lambda}
\end{equation*}
such that
    $X_h' \mapsto X_h$
for every $h \not\in \{i,j,k,l,m\}$, $X_i' \mapsto X_i^{-1}$ and:
\begin{description}
\item[Case 1] if the edges $\lambda_j$, $\lambda_k$, $\lambda_l$,
$\lambda_m$  are distinct, then
\begin{align*}
X'_j &\mapsto (1+qX_i)X_j \qquad\, X'_k \mapsto
(1+qX_i^{-1})^{-1}X_k \\
    X'_l &\mapsto (1+qX_i)X_l
\qquad
    X'_m
\mapsto (1+qX_i^{-1})^{-1}X_m ;
\end{align*}

\item[Case 2] if $\lambda_j$ is identified with $\lambda_k$, and
$\lambda_l$ is distinct from $\lambda_m$, then
$$ X'_j  \mapsto   X_iX_j \quad X'_l
\mapsto (1+qX_i)X_l \quad X'_m \mapsto (1+qX_i^{-1})^{-1}X_m
$$

\item[Case 3] (the inverse of Case 2) if $\lambda_j$ is identified
with $\lambda_m$, and $\lambda_k$ is distinct from $\lambda_l$,
then
$$X_j\mapsto   X_iX_j \quad X'_k
\mapsto (1+qX_i^{-1})^{-1}X_k \quad X'_l  \mapsto (1+qX_i)X_l
$$

\item[Case 4] if $\lambda_j$ is identified with $\lambda_l$, and
$\lambda_k$ is distinct from $\lambda_m$, then
\begin{gather*}
X'_j  \mapsto (1+qX_i) (1+q^3X_i) X_j \\ X'_k \mapsto
(1+qX_i^{-1})^{-1}X_k \quad X'_m \mapsto (1+qX_i^{-1})^{-1}X_m
\end{gather*}

\item[Case 5] (the inverse of Case 4) if $\lambda_k$ is identified
with $\lambda_m$, and $\lambda_j$ is distinct from $\lambda_l$,
then
\begin{gather*}
    X'_j  \mapsto  (1+qX_i)X_j
\quad X'_l
    \mapsto  (1+qX_i)X_l \\
    X'_k \mapsto
(1+qX_i^{-1})^{-1}(1+q^3X_i^{-1})^{-1}X_k
\end{gather*}

\item[Case 6] if $\lambda_j$ is identified with $\lambda_k$, and
$\lambda_l$ is identified with $\lambda_m$ (in which case $S$ is a
$3$--times punctured sphere), then
$$
    X'_j  \mapsto  X_iX_j \quad X'_l  \mapsto
X_iX_l;
$$

\item[Case 7] (the inverse of Case 6) if $\lambda_j$ is identified
with $\lambda_m$, and $\lambda_k$ is identified with $\lambda_l$
(in which case $S$ is a $3$--times punctured sphere), then
$$
    X'_j  \mapsto  X_iX_j
\quad X'_k  \mapsto  X_iX_k;
$$

\item[Case 8] if $\lambda_j$ is identified with $\lambda_l$, and
$\lambda_k$ is identified with $\lambda_m$ (in which case $S$ is a
once punctured torus), then
\begin{gather*}
    X'_j  \mapsto
(1+qX_i)(1+q^3X_i)X_j \\
    X'_k
\mapsto (1+qX_i^{-1})^{-1}(1+q^3X_i^{-1})^{-1}X_k
\end{gather*}

\end{description}
\end{prop}

\begin{proof} By inspection, these formulas
are compatible with the skew-commutativity relations
(\ref{eqn:SkewCom}). For instance, in Case~4, $X_i'X_j' = q^4
X_j'X_i'$, $X_i'X_k' = q^{-2} X_k'X_i'$ and $X_j'X_k = q^2
X_k'X_j'$ in ${\mathcal{T}}^q_{\lambda'}$, and
\begin{gather*}
\bigl [ X_i^{-1} \bigr] \bigl [ (1+qX_i)(1+q^3X_i) X_j \bigr] =
q^4\bigl [ (1+qX_i)(1+q^3X_i) X_j \bigr]
\bigl [ X_i^{-1} \bigr]\\
\bigl [ X_i^{-1} \bigr] \bigl [ (1+qX_i^{-1})^{-1}X_k \bigr] =
q^{-2}\bigl [ (1+qX_i^{-1})^{-1}X_k \bigr]
\bigl [ X_i^{-1} \bigr]\\
\bigl [ (1+qX_i)(1+q^3X_i) X_j \bigr] \bigl [(1+qX_i^{-1})^{-1}X_k
\bigr] \qquad\qquad\qquad\qquad\qquad\qquad
\\
\qquad\qquad\qquad \qquad\qquad\qquad =q^2\bigl [
(1+qX_i^{-1})^{-1}X_k \bigr] \bigl [ (1+qX_i)(1+q^3X_i) X_j \bigr]
\end{gather*}
in $\widehat{\mathcal{T}}^q_{\lambda}$ since $X_iX_j = q^{-4}
X_jX_i$, $X_iX_k = q^2X_kX_j$ and $X_jX_k = X_kX_j$.  It follows
that these formulas extend to a unique algebra homomorphism
$\phi_{\lambda\lambda'}^q: {\mathcal{T}}^q_{\lambda'} \rightarrow
\widehat{\mathcal{T}}^q_{\lambda}$.

Extending $\phi_{\lambda\lambda'}^q$ to the fraction division
algebra $\widehat{\mathcal{T}}^q_{\lambda'}$ will require a little
care. Indeed, for a formal fraction $PQ^{-1}\in
\widehat{\mathcal{T}}^q_{\lambda'}$ with $P$, $Q\in
{\mathcal{T}}^q_{\lambda'}$, we want to define
$\Phi_{\lambda\lambda'}^q(PQ^{-1})= \phi_{\lambda\lambda'}^q (P)
\phi_{\lambda\lambda'}^q(Q)^{-1}\in
\widehat{\mathcal{T}}^q_{\lambda}$. For this we need
$\phi_{\lambda\lambda'}^q(Q)$ to be non-zero. In other words, we
need to show that $\phi_{\lambda\lambda'}^q$ is injective.

For this, we take advantage of the reflexivity of the above
formulas. Exchanging the r\^{o}les of $\lambda$ and $\lambda'$, we
also have an algebra homomorphism $\phi_{\lambda'\lambda}^q:
{\mathcal{T}}^q_{\lambda} \rightarrow
\widehat{\mathcal{T}}^q_{\lambda'}$. Note that the elements of
$\phi_{\lambda\lambda'}^q \bigl( {\mathcal{T}}^q_{\lambda'}
\bigr)$ are all polynomials in the $X_h^{\pm1}$,with $h=1$, 2,
\dots, $n$ and in the quantity $(1+qX_i^{-1})^{-1}$. After
checking the skew-commutativity relations and because
$\phi_{\lambda'\lambda}^q (X_i) = (X_i')^{-1}$, we can therefore
extend $\phi_{\lambda'\lambda}^q$ to $\phi_{\lambda\lambda'}^q
\bigl( {\mathcal{T}}^q_{\lambda'} \bigr)$ by defining
$\phi_{\lambda'\lambda}^q \bigl( (1+qX_i^{-1})^{-1} \bigr)
=(1+qX_i')^{-1}$.

We now have a composition
$$\mathcal{T}^q_{\lambda'}\
\stackrel{\phi_{\lambda\lambda'}^q}{ \longrightarrow}\
\phi_{\lambda\lambda'}^q \bigl( {\mathcal{T}}^q_{\lambda'} \bigr)\
\stackrel{\phi_{\lambda'\lambda}^q} {\longrightarrow}\
\widehat{\mathcal{T}}^q_{\lambda'}$$ of two algebra homomorphisms.
By inspection of the formulas defining these homomorphisms,
$\phi_{\lambda'\lambda}^q \circ \phi_{\lambda\lambda'}^q (X_h') =
X_h'$ for every generator $X_h'$ of $\mathcal{T}^q_{\lambda'}$. It
follows that $\phi_{\lambda'\lambda}^q \circ
\phi_{\lambda\lambda'}^q$ is the inclusion map. In particular, it
is injective, which proves that the first homomorphism
$\phi_{\lambda\lambda'}^q: {\mathcal{T}}^q_{\lambda'} \rightarrow
\widehat{\mathcal{T}}^q_{\lambda}$ is injective.

As indicated before, this enables us to extend
$\phi_{\lambda\lambda'}^q: {\mathcal{T}}^q_{\lambda'} \rightarrow
\widehat{\mathcal{T}}^q_{\lambda}$ to a map
$\Phi_{\lambda\lambda'}^q: \widehat{\mathcal{T}}^q_{\lambda'}
\rightarrow \widehat{\mathcal{T}}^q_{\lambda}$ by setting
$\Phi_{\lambda\lambda'}^q(PQ^{-1})= \phi_{\lambda\lambda'}^q (P)
\phi_{\lambda\lambda'}^q(Q)^{-1}$ if $P$, $Q\in
{\mathcal{T}}^q_{\lambda'}$. It can be shown that this map
$\Phi_{\lambda\lambda'}^q$ is an algebra homomorphism. See
\cite[Chap.~9]{GooWar}.

A symmetric argument provides an algebra homomorphism
$\Phi_{\lambda'\lambda}^q: \widehat{\mathcal{T}}^q_{\lambda}
\rightarrow \widehat{\mathcal{T}}^q_{\lambda'}$. From the fact
that $\phi_{\lambda'\lambda}^q \circ \phi_{\lambda\lambda'}^q
(X_h') = X_h'$ for every generator $X_h'$ of
$\mathcal{T}^q_{\lambda'}$, we conclude that
$\Phi_{\lambda'\lambda}^q \circ \Phi_{\lambda\lambda'}^q$ is the
identity. By symmetry, it follows that $\Phi_{\lambda\lambda'}^q:
\widehat{\mathcal{T}}^q_{\lambda'} \rightarrow
\widehat{\mathcal{T}}^q_{\lambda}$ is an isomorphism, with inverse
$\Phi_{\lambda'\lambda}^q$.
\end{proof}

When $\lambda'$ is obtained from $\lambda$ by an edge reindexing,
we also have the immediate relation.

\begin{prop}
\label{prop:QuantumReindex} Suppose that the ideal triangulations
$\lambda$, $\lambda'\in \Lambda(S)$ are obtained from each other
by an edge reindexing, namely that $\lambda'_i =
\lambda_{\alpha(i)}$ for some permutation $\alpha \in \mathfrak
S_n$. Then there exists a unique isomorphism
$$\Phi_{\lambda\lambda'}^q:
\widehat{\mathcal{T}}^q_{\lambda'} \rightarrow
\widehat{\mathcal{T}}^q_{\lambda}$$ such that
$\Phi_{\lambda\lambda'}^q(X_i') = X_{\alpha(i)}$ for every
generator $X_i'$ of $\widehat{\mathcal{T}}^q_{\lambda'}$.  \qed
\end{prop}

The following facts are straightforward from the definitions.
\begin{prop}
\label{prop:QuantumRelations} The isomorphisms
$\Phi_{\lambda\lambda'}^q: \widehat{\mathcal{T}}^q_{\lambda'}
\rightarrow \widehat{\mathcal{T}}^q_{\lambda}$ defined by
Propositions~\ref{prop:QuantumDiagExch} and
\ref{prop:QuantumReindex} satisfy the following properties.
\begin{enumerate}

\item Composition Relation: $\Phi_{\lambda\beta(\lambda)}^q \circ
\Phi_{\beta(\lambda) \,\alpha\circ\beta(\lambda)} ^q=
\Phi_{\lambda \,\alpha\circ\beta(\lambda)}^q$ for every $\alpha$,
$\beta\in \mathfrak S_n$.

\item Reflexivity Relation: $\Phi_{\lambda\Delta_i(\lambda)}^q
\circ \Phi_{\Delta_i(\lambda)\lambda}^q=\mathrm{Id}$

\item Reindexing Relation: $\Phi_{\lambda\alpha(\lambda)}^q \circ
\Phi_{\alpha(\lambda) \,\Delta_i\circ\alpha(\lambda)} ^q
=\Phi_{\lambda\Delta_{\alpha(i)}(\lambda)}^q \circ
\Phi_{\Delta_{\alpha(i)}(\lambda)
\,\alpha\circ\Delta_{\alpha(i)}(\lambda)} ^q$ for every $\alpha\in
\mathfrak S_n$.

\item Distant Commutativity Relation:
$\Phi_{\lambda\Delta_j(\lambda)}^q \circ \Phi_{\Delta_j(\lambda)
\,\Delta_i\circ\Delta_j(\lambda)} ^q=
\Phi_{\lambda\Delta_i(\lambda)}^q \circ \Phi_{\Delta_i(\lambda)
\,\Delta_j\circ\Delta_i(\lambda)} ^q$ if the edges $\lambda_i$,
$\lambda_j$ of $\lambda$ do not belong to a same triangle. \qed

\end{enumerate}
\end{prop}

Dealing with the Pentagon Relation will require more efforts.

\section{The Pentagon Relation}
\label{sect:Pentagon}

The goal of this section is to show that the isomorphisms
$\Phi_{\lambda\lambda'}^q$ constructed in the previous section are
compatible with the pentagon relation satisfied by the diagonal
exchanges $\Delta_i$.

The main tool to achieve this is borrowed from \cite{FocChe}. It
consists of a certain order 5 automorphism of the quantum torus.
The quantum torus $\mathcal W^q$ is the algebra defined by the
generators $U^{\pm1}$, $V^{\pm1}$ and by the relation $VU=q^2UV$.

Set $U_{(0)} = U$, $V_{(0)} = V$ and inductively define
\begin{equation}
\label{eqn:cycle}
    U_{(k+1)}  = (1+qU_{(k)})V_{(k)},
\quad V_{(k+1)}  = U_{(k)}^{-1}
\end{equation}
   From this definition we
deduce the following relations.
\begin{gather}
\label{eqn:commutation} U_{(k+1)}U_{(k)}=q^2U_{(k)}U_{(k+1)},\quad
V_{(k)}U_{(k)}=q^2U_{(k)}V_{(k)}\\
\label{eqn:1+qU} U_{(k+1)}U_{(k-1)}=1+qU_{(k)},\quad
U_{(k-1)}U_{(k+1)}=1+q^{-1}U_{(k)}
\end{gather}
This provides us with a few combinatorial tricks to deal with
noncommutativity. One of them is that the products
$U_{(k+1)}U_{(k-1)}$ and $U_{(k-1)}U_{(k+1)}$ both commute with
$U_{(k)}$. Also, the two equations of (\ref{eqn:1+qU}) deform to:
\begin{gather}
\label{decomp1} q^{-1}U_{(k)}^{-1}U_{(k+1)}U_{(k-1)}
= 1+q^{-1}U_{(k)}^{-1}\\
\label{decomp2} qU_{(k)}^{-1}U_{(k-1)} U_{(k+1)}=1+ qU_{(k)}^{-1}.
\end{gather}
What we have accomplished is decomposing the polynomial factors
$1+q^{\pm1}U_{(k)}^{\pm1}$ into monomials. This decomposition will
turn out to be very convenient later on.

\begin{lem}[Chekhov-Fock]
\label{lem:order5} $U_{(k+5)}=U_{(k)}$ for every $k\in\mathbb{Z}$.
\end{lem}

\begin{proof} This immediately follows from
iterated applications of (\ref{eqn:cycle}):
\begin{align*} U_{(k+5)} & =
(1+qU_{(k+4)})U_{(k+3)}^{-1}\\
    & =
U_{(k+3)}^{-1} +q\left((1+qU_{(k+3)})U_{(k+2)}^{-1}\right)
U_{(k+3)}^{-1}
\\
& = U_{(k+3)}^{-1}+qU_{(k+2)}^{-1}U_{(k+3)}^{-1}
    +q^2U_{(k+3)}U_{(k+2)}^{-1}
U_{(k+3)}^{-1}\\
& = U_{(k+3)}^{-1} +qU_{(k+2)}^{-1}U_{(k+3)}^{-1}
    +U_{(k+2)}^{-1}
\end{align*}
and
\begin{align*}
    U_{(k)} & =
U_{(k+2)}^{-1}(1+qU_{(k+1)})\\
& = U_{(k+2)}^{-1} +qU_{(k+2)}^{-1}
\left(U_{(k+3)}^{-1}(1+qU_{(k+2)})\right)\\
& = U_{(k+2)}^{-1}+qU_{(k+2)}^{-1}U_{(k+3)}^{-1}
+q^2U_{(k+2)}^{-1}U_{(k+3)}^{-1}U_{(k+2)}\\
    & =
U_{(k+2)}^{-1} +qU_{(k+2)}^{-1}U_{(k+3)}^{-1}+U_{(k+3)}^{-1}
\end{align*}
\vskip -\belowdisplayskip \vskip -\baselineskip
\end{proof}

Consider a pentagon cycle of geodesic laminations $\lambda_{(0)}$,
$\lambda_{(1)} = \Delta_i(\lambda_{(0)})$, $\lambda_{(2)} =
\Delta_j(\lambda_{(1)})$, $\lambda_{(3)} =
\Delta_i(\lambda_{(2)})$, $\lambda_{(4)} =
\Delta_j(\lambda_{(3)})$, $\lambda_{(5)} = \Delta_i(\lambda_{(4)})
=\alpha_{ i \leftrightarrow j}(\lambda_{(0)})$ as in
Figure~\ref{fig:pentagon}. We want to show that
$$
\Phi_{\lambda(0)\lambda(1)}^q \circ \Phi_{\lambda(1)\lambda(2)}^q
\circ \Phi_{\lambda(2)\lambda(3)}^q \circ
\Phi_{\lambda(3)\lambda(4)}^q \circ \Phi_{\lambda(4)\lambda(5)}^q
= \Phi_{\lambda(0)\lambda(5)}^q
$$
where the $\Phi_{\lambda(k)\lambda(k+1)}^q$ are defined by
Proposition~\ref{prop:QuantumDiagExch}, and
$\Phi_{\lambda(0)\lambda(5)}^q$ by
Proposition~\ref{prop:QuantumReindex}.

Let $\lambda_a$, $\lambda_b$, $\lambda_c$, $\lambda_d$,
$\lambda_e$ denote the sides of the pentagon, and let $\lambda_i$
and $\lambda_j$ be the diagonal edges of $\lambda_{(0)}$, all
labelled as in Figure~\ref{fig:pentagon}.

Note that the $\lambda_a$, $\lambda_b$, $\lambda_c$, $\lambda_d$,
$\lambda_e$ are edges of all the $\lambda_{(k)}$, and in
particular are associated to generators $X_a$, $X_b$, $X_c$,
$X_d$, $X_e $ of $\mathcal T^q_{\lambda_{(i)}}$. Let $A_{(k)}$,
$B_{(k)}$, $C_{(k)}$, $D_{(k)}$, $E_{(k)} \in \widehat{\mathcal
T}^q_{\lambda_{(0)}}$ denote the respective images of $X_a$,
$X_b$, $X_c$, $X_d$, $X_e \in \mathcal T^q_{\lambda_{(k)}}$ under
$\Phi^q_{\lambda_{(0)}\lambda_{(1)}} \circ \dots \circ
\Phi^q_{\lambda_{(k-1)}\lambda_{(k)}}$.

We do the same thing with diagonals, but with an additional twist.
The diagonals of the pentagon in $\lambda_{(k)}$ are always its
$i$--th and $j$--th edges, but their relative configuration with
respect to each other alternates. In particular, the factor
$\sigma_{ij}^{\lambda_{(k)}}$ occurring in the skew-commutativity
relations (\ref{eqn:SkewCom}) is equal to $(-1)^{k+1}$. For this
reason, we let $U_{(k)}$ and $V_{(k)} \in \widehat{\mathcal
T}^q_{\lambda_{(0)}}$ be the respective images of $X_{\alpha^k_{i
\leftrightarrow j}(i)}$  and $X_{\alpha^k_{i \leftrightarrow
j}(j)}\in \mathcal T^q_{\lambda_{(i)}}$ under
$\Phi^q_{\lambda_{(0)}\lambda_{(1)}} \circ \dots \circ
\Phi^q_{\lambda_{(k-1)}\lambda_{(k)}}$, where $\alpha^k_{i
\leftrightarrow j} \in \mathfrak S_n$ is the transposition
exchanging $i$ and $j$.  This is specially designed so that
$V_{(k)}U_{(k)} = q^2 U_{(k)}V_{(k)}$ for every $k$. In addition,
as one moves from $\lambda_{(k)}$ to $\lambda_{(k+1)}$, the
diagonal exchange is always performed on the edge corresponding to
$U_{(k)}$.

   From the definition of
$\Phi^q_{\lambda_{(k-1)}\lambda_{(k)}}$, the elements $U_{(k)}$,
$V_{(k)}$ satisfy the induction relation (\ref{eqn:cycle}). In
particular, it follows from Lemma~\ref{lem:order5} that
$U_{(k+5)}=U_{(k)}$ and $V_{(k+5)}=V_{(k)}$.

The induction formulas for the other elements $A_{(k)}$,
$B_{(k)}$, $C_{(k)}$, $D_{(k)}$, $E_{(k)}$ depend on whether there
are identifications between the sides of the pentagon. We first
consider the case of an \emph{embedded} polygon, with no
identification between its sides.

\begin{prop}[Chekhov-Fock]
\label{prop:chefock} The Pentagon Relation
$$
\Phi_{\lambda(0)\lambda(1)}^q \circ \Phi_{\lambda(1)\lambda(2)}^q
\circ \Phi_{\lambda(2)\lambda(3)}^q \circ
\Phi_{\lambda(3)\lambda(4)}^q \circ \Phi_{\lambda(4)\lambda(5)}^q
= \Phi_{\lambda(0)\lambda(5)}^q
$$
is satisfied in the case of an embedded pentagon.
\end{prop}

\begin{proof} Set
$ \Psi= \Phi_{\lambda(0)\lambda(1)}^q \circ
\Phi_{\lambda(1)\lambda(2)}^q \circ \Phi_{\lambda(2)\lambda(3)}^q
\circ \Phi_{\lambda(3)\lambda(4)}^q \circ
\Phi_{\lambda(4)\lambda(5)}^q $ to simplify the notation. We need
to show that $\Psi(X_k) = X_{\alpha_{i \leftrightarrow j}(k)}$ for
every $k$. The property is immediate for those $X_k$ which
correspond to edges outside of the pentagon.

By definition of the $U_{(k)}$, $V_{(k)}$ and by
Lemma~\ref{lem:order5}, $\Psi(X_i)= V_{(5)} = V_{(0)}=X_j$ and
$\Psi(X_j)= U_{(5)} = U_{(0)}=X_i$. This proves the property for
$X_i$ and $X_j$.

For the sides of the pentagon, $\Psi(X_a)= A_{(5)}$, $\Psi(X_b)=
B_{(5)}$, $\Psi(X_c)= C_{(5)}$, $\Psi(X_d)= D_{(5)}$ and
$\Psi(X_e)= E_{(5)}$.
    From the definition of
$\Phi^q_{\lambda_{(i-1)}\lambda_{(i)}}$ in
Proposition~\ref{prop:QuantumDiagExch},
\begin{align*}
A_{(5)} & = \big(1+qU_{(3)}^{-1}\big)^{-1}
\big(1+qU_{(2)}^{-1}\big)^{-1} \big(1+qU_{(0)}\big)A_{(0)}
\\
    B_{(5)} & =
\big(1+qU_{(4)}^{-1}\big)^{-1} \big(1+qU_{(2)}\big)
\big(1+qU_{(0)}^{-1}\big)^{-1}B_{(0)}
\\
C_{(5)} & = \big(1+qU_{(4)}\big) \big(1+qU_{(2)}^{-1}\big)^{-1}
\big(1+qU_{(1)}^{-1}\big)^{-1}C_{(0)}
    \\
D_{(5)} & = \big(1+qU_{(4)}^{-1}\big)^{-1}
\big(1+qU_{(3)}^{-1}\big)^{-1}
\big(1+qU_{(1)}\big)D_{(0)}\\
    E_{(5)} & =
\big(1+qU_{(3)}\big)\big(1+qU_{(1)}^{-1}
\big)^{-1}\big(1+qU_{(0)}^{-1}\big)^{-1} E_{(0)}
\end{align*}

Using (\ref{decomp2}),
\begin{align}
\big(1+qU_{(k+3)}^{-1}\big)^{-1} \big(1+qU_{(k+2)}^{-1}\big)^{-1}
    & =  q^{-1}U_{(k+3)}U_{(k+4)}^{-1}
U_{(k+2)}^{-1}q^{-1}U_{(k+2)}U_{(k+3)}^{-1}
U_{(k+1)}^{-1}\nonumber
\\
    & =
    U_{(k+4)}^{-1}U_{(k+1)}^{-1}\nonumber
\\
& = U_{(k-1)}^{-1}U_{(k+1)}^{-1}\nonumber
\\
    &
=  (1+qU_{(k)})^{-1}. \label{eqn:ThreeTerm}
\end{align}
Combining this relation with $U_{(k'+5)}=U_{(k')}$, we get that
$A_{(5)}=A_{(0)}=X_a$, $B_{(5)}=B_{(0)}=X_b$,
$C_{(5)}=C_{(0)}=X_c$, $D_{(5)}=D_{(0)}=X_d$ and
$E_{(5)}=E_{(0)}=X_e$ as required.
\end{proof}

We now have to worry about possible identifications between the
sides of the pentagon. Note that it suffices to prove the property
for any cyclic permutation of the $\lambda_{(k)}$. This reduces
the analysis to 6 possible cases.
\begin{enumerate}
\item no identification (embedded pentagon); \item $\lambda_a =
\lambda_b$, and the other sides are distinct; \item $\lambda_a =
\lambda_c$, and the other sides are distinct; \item $\lambda_a =
\lambda_b$ and $\lambda_c = \lambda_d$; \item $\lambda_a =
\lambda_b$ and $\lambda_c = \lambda_e$; \item $\lambda_a =
\lambda_c$ and $\lambda_b = \lambda_e$.
\end{enumerate}

We already considered Case~1.

In the other cases, note that the identifications between the
sides of the pentagon have no impact on the images of $X_i$,
$X_j$. We therefore only need to consider $X_a$, $X_b$, $X_c$,
$X_d$, $X_e$.

\medskip
\noindent\textbf{Case 2}: $\lambda_a = \lambda_b$, and the other
sides are distinct.

Using (\ref{eqn:ThreeTerm}) and (\ref{eqn:1+qU}),
\begin{align*}
\Psi(X_a)= A_{(5)} & =   \big(1+qU_{(4)}^{-1}\big)^{-1}
\big(1+qU_{(3)}^{-1}\big)^{-1} U_{(2)}U_{(0)}A_{(0)}
\\
& = (1+qU_{(1)})^{-1}U_{(2)}U_{(0)}A_{(0)}
\\
& =  A_{(0)} =X_a
\end{align*}
The argument for $X_c$, $X_d$, $X_e$ is identical to that of
Proposition~\ref{prop:chefock}.

\medskip
\noindent\textbf{Case 3}:  $\lambda_a = \lambda_c$, and the other
sides are distinct.

Using (\ref{eqn:ThreeTerm}) and (\ref{eqn:commutation}),
\begin{align*}
A_{(5)} & = \big(1+qU_{(4)}\big)\big(1+qU_{(3)}^{-1} \big)^{-1}
\big(1+qU_{(2)}^{-1}\big)^{-1} \big(
1+q^3U_{(2)}^{-1}\big)^{-1}A_{(2)}
\\
    & =
\big(1+qU_{(4)}\big) \big(1+qU_{(0)}\big)^{-1} \big(
1+q^3U_{(2)}^{-1}\big)^{-1}A_{(2)}
\\
    & =
U_{(0)}U_{(3)}U_{(4)}^{-1}U_{(1)}^{-1}\big(
1+q^3U_{(2)}^{-1}\big)^{-1}A_{(2)}
\\
&= U_{(0)}U_{(3)}U_{(4)}^{-1} \big(1+qU_{(2)}^{-1}\big)^{-1}
U_{(1)}^{-1}A_{(2)}
\\
    & =
U_{(0)}U_{(3)}U_{(4)}^{-1}\,q^{-1} U_{(3)}^{-1}U_{(1)}^{-1}U_{(2)}
U_{(1)}^{-1}A_{(2)}
\end{align*}
and
\begin{align*} A_{(2)}
& = \big(1+qU_{(1)}^{-1}\big)^{-1} \big(1+qU_{(0)}\big)A_{(0)}
\\
    & =
q^{-1}U_{(2)}^{-1}U_{(0)}^{-1}U_{(1)} U_{(1)}U_{(4)}A_{(0)}
\end{align*}
Then $$
A_{(5)}A_{(0)}^{-1}=q^{-2}U_{(0)}U_{(3)}U_{(4)}^{-1}U_{(3)}^{-1}
U_{(1)}^{-1}U_{(2)}U_{(1)}^{-1}U_{(2)}^{-1}U_{(0)}^{-1}U_{(1)}
U_{(1)}U_{(4)}=\mathbf{1}.$$ It follows that $\Psi(X_a) = A_{(5)}
= A_{(0)} = X_a$. The argument for $X_b$, $X_d$, $X_e$ is
identical to that for Case~1.

\medskip
\noindent\textbf{Case 4}:
    $\lambda_a = \lambda_b$ and
$\lambda_c = \lambda_d$.

The fact that $\Psi(X_a)=X_a$ is proved as in Case~2. Shifting
indices by 2 and rotating the picture, the argument of Case~2
gives $C_{(7)}=C_{(2)}$. We can then backtrack to
    $C_{(5)}=C_{(0)}$ by using the 5--periodicity of
the $U_{(k)}$. This proves that $\Psi(X_c)=X_c$. The remaining
generator $X_e$ is treated as in Case~1.

\medskip
\noindent\textbf{Case 5}: $\lambda_a = \lambda_b$ and $\lambda_c =
\lambda_e$.

Again $\Psi(X_a)=X_a$ as in Case~2, and $\Psi(X_d)=X_d$ as in
Case~1. The fact that $\Psi(X_c)=X_c$ follows from Case~3 after
shifting indices by 2 as above.

\medskip
\noindent\textbf{Case 6}: $\lambda_a = \lambda_c$ and $\lambda_b =
\lambda_e$. This again follows from a combination of the arguments
of Cases~1, 2 and 3.
\medskip

In conclusion, we have proved:
\begin{prop}
\label{prop:pentagon} The Pentagon Relation
$$
\Phi_{\lambda(0)\lambda(1)}^q \circ \Phi_{\lambda(1)\lambda(2)}^q
\circ \Phi_{\lambda(2)\lambda(3)}^q \circ
\Phi_{\lambda(3)\lambda(4)}^q \circ \Phi_{\lambda(4)\lambda(5)}^q
= \Phi_{\lambda(0)\lambda(5)}^q
$$
is satisfied in all cases. \qed
\end{prop}

\section{The quantum Teichm\"uller space}
\label{sect:QuantumTeich}

We can now state and prove the main result of this paper.

\begin{thm}
\label{thm:main} There is a unique family of algebra isomorphisms
$$\Phi_{\lambda\lambda'}^q:
\widehat{\mathcal{T}}^q_{\lambda'} \rightarrow
\widehat{\mathcal{T}}^q_{\lambda}$$ defined as $\lambda$,
$\lambda' \in \Lambda(S)$ ranges over all pairs of ideal
triangulations, such that:
\begin{enumerate}
\item $\Phi_{\lambda\lambda''}^q = \Phi_{\lambda\lambda'}^q \circ
\Phi_{\lambda'\lambda''}^q$ for every $\lambda$, $\lambda'$,
$\lambda''\in \Lambda(S)$;

\item $\Phi_{\lambda\lambda'}^q$ is the isomorphism of
Proposition~\ref{prop:QuantumDiagExch} when $\lambda'$ is obtained
from $\lambda$ by a diagonal exchange;

\item $\Phi_{\lambda\lambda'}^q$ is the isomorphism of
Proposition~\ref{prop:QuantumReindex} when $\lambda'$ is obtained
from $\lambda$ by an edge reindexing.
\end{enumerate}
\end{thm}
\begin{proof} Use Theorem~\ref{thm:Penner1} to
connect $\lambda$ to $\lambda'$ by a sequence
$\lambda=\lambda_{(0)}$, $\lambda_{(1)}$, \dots,
$\lambda_{(m)}=\lambda'$ where each $\lambda_{(k+1)}$ is obtained
from $\lambda_{(k)}$ by a diagonal exchange or by an edge
reindexing, and define $\Phi_{\lambda\lambda'}^q$ as the
composition of the $\Phi_{\lambda_{(k)}\lambda_{(k+1)}}^q$
provided by Propositions~\ref{prop:QuantumDiagExch} and
\ref{prop:QuantumReindex}. Theorem~\ref{thm:Penner2} and
Propositions~\ref{prop:QuantumRelations} and \ref{prop:pentagon}
show that this $\Phi_{\lambda_{(k)}\lambda_{(k+1)}}^q$ is
independent of the choice of the sequence of $\lambda_{(k)}$.

The uniqueness immediately follows from Theorem~\ref{thm:Penner1}.
\end{proof}

The \emph{quantum (enhanced) Teichm\"uller space} of $S$ can now
be defined as the algebra
$$
\widehat {\mathcal{T}}^q_S= \bigg(
\bigsqcup_{\lambda\in\Lambda(S)}
\widehat{\mathcal{T}}^q_{\lambda}\bigg)/\sim
$$
where the relation $\sim$ is defined by the property that, for
$X\in \widehat{\mathcal{T}}^q_{\lambda}$ and $X'\in
\widehat{\mathcal{T}}^q_{\lambda'}$,
$$
X \sim X' \Leftrightarrow X=\Phi^q_{\lambda,\lambda'}(X').
$$

Note that the definition is specially designed so that, when
$q=1$, there is a natural isomorphism between
$\widehat{\mathcal{T}}^1_S$ and the algebra
$\mathrm{Rat}\,\widetilde{\mathcal T}(S)$ of rational functions on
the enhanced Teichm\"uller space $\widetilde{\mathcal T}(S)$.

\section{The quantum cusped Teichm\"uller space}
\label{sect:QuantumCusped}

The \emph{cusped Teichm\"uller space} $\mathcal {CT}(S)$ is the
set of isotopy classes of hyperbolic metrics on $S$ for which all
ends are of cusp type. Since such a metric admits a unique
enhancement, $\mathcal {CT}(S)$ is a natural subspace of the
enhanced Teichm\"uller space $\widetilde{\mathcal T}(S)$.

Consider the Thurston parametrization $\phi_\lambda
:\widetilde{\mathcal T}(S) \rightarrow \mathbb R_+^n$ of
$\widetilde{\mathcal T}(S)$ by the shear coordinates $x_i$
associated to the edges $\lambda_i$ of an ideal triangulation
$\lambda \in \Lambda(S)$. Recall that we associated to the $j$--th
puncture $v_j$ of $S$ the exponential length parameter
$p_j=x_1^{k_{1j}}x_2^{k_{2j}}\dots x_n^{k_{nj}}$ where $k_{ij}\in
\{0,1,2\}$ is the number of end points of the edge $\lambda_i$
that are equal to $v_j$. We observed in \S~\ref{sect:Shear} that
$\mathcal {CT}(S) \subset \widetilde{\mathcal T}(S)$ corresponds
under $\phi_\lambda$ to the set of those $x\in  \mathbb R_+^n$
such that $p_j=1$ for every $j=1$, \dots, $p$.

It is also natural to consider the  product $h=x_1 x_2 \dots x_n$.
Note that $h^2 = p_1 p_2 \dots p_p$ since every edge $\lambda_i$
has two end points (so that $\sum_j k_{ij}=2$ for every $i$). In
particular, $h$ is identically 1 on the image of $\mathcal
{CT}(S)$ under $\phi_\lambda$ since the $x_i$ are all real
positive.

In this non-quantum context, consider the elements
$P_j=X_1^{k_{1j}}X_2^{k_{2j}}\dots X_n^{k_{nj}}$ and $H=X_1 X_2
\dots X_n \in \mathbb C (X_1, X_2, \dots, X_n)$. Note that again
$H^2 = P_1 P_2 \dots P_p$.

It can be shown that there exists an isomorphism
$$
\mathbb C (X_1, X_2, \dots, X_n) \rightarrow \mathbb C (Y_1, Y_2,
\dots, Y_n)$$
    sending the first $p-1$ elements $P_1$,
$P_2$, \dots, $P_{p-1}$ to $Y_1$, $Y_2$, \dots, $Y_{p-1}$ and $H$
to $Y_p$, respectively. See for instance \cite[\S 3]{BonLiu}. It
follows that the isomorphism
$$\Phi_\lambda :
\mathrm{Rat}\,\widetilde{\mathcal T}(S) \rightarrow\mathbb C(X_1,
X_2, \dots, X_n)$$ defined by $\phi_\lambda$ induces an
isomorphism
$$\Psi_\lambda :
\mathrm{Rat}\,\mathcal {CT}(S) \rightarrow\mathbb C(X_1, X_2,
\dots, X_n)/I$$ where $I$ is the ideal generated by the elements
$P_j-1$ and $H-1$.  Note that it is important to include $H-1$ to
make sure that $I$ is prime.

We extend this to the quantum set-up by introducing the elements
$$P_j=q^{-\sum_{i<i'}
\sigma_{ii'}^\lambda k_{ij}k_{i'j}} X_1^{k_{1j}}X_2^{k_{2j}}\dots
X_n^{k_{nj}}$$ and
$$H= q^{-\sum_{i<i'}
\sigma_{ii'}^\lambda } X_1 X_2 \dots X_n$$ of $\mathcal
T_\lambda^q$. The $q$--factor is introduced for the following
property.

\begin{prop}
\label{prop:InvarReindex} The elements $H$ and $P_j\in \mathcal
T_\lambda^q$ are each invariant under reindexing of the edges of
$\lambda$.
\end{prop}

\begin{proof}
It suffices to verify this for the transposition
$\alpha_{k\leftrightarrow k+1}$. Exchanging $k$ and $k+1$ in the
product $X_1^{k_{1j}}X_2^{k_{2j}}\dots X_n^{k_{nj}}$ of $P_j$
results in multiplying this element by
$q^{-2\sigma_{k,k+1}^\lambda k_{kj}k_{k+1\,j}}$. However, in the
exponent ${-\sum_{i<i'} \sigma_{ii'}^\lambda k_{ij}k_{i'j}}$ of
the $q$--factor, $\sigma_{k\,k+1}^\lambda$ is replaced by
$\sigma_{k+1\,k}^\lambda =-\sigma_{k\,k+1}^\lambda$, so that this
exponent increases by $2\sigma_{k,k+1}^\lambda$ $k_{kj}
k_{k+1\,j}$. Consequently, the two contributions cancel out and
$P_j$ remains invariant. The argument is similar for $H$.
\end{proof}

\begin{rem}
The $q$--factor in the definition of $H$ and the $P_j$ is
traditionally known as the quantum ordering in the physics
literature. It is better explained in a situation when we can
write $X_i= \exp \xi_i$,  where the $\xi_i$ are symbols with
central commutators $\xi_i\xi_j-\xi_j\xi_i=2\pi \mathrm
i\hbar\sigma_{ij} \in \mathbb C$. It then follows from the
Campbell-Hausdorff formula that
\begin{align*}
\exp(k_1\xi_1+k_2\xi_2+ \cdots +k_n \xi_n) &=\mathrm e^{- \pi
\mathrm i \hbar \sum_{i<i'} \sigma_{ii'}k_i k_{i'}} \exp(k_1\xi_1)
\exp(k_2\xi_2)\dots \exp(k_n\xi_n)
\\
&= q^{-\sum_{i<i'} \sigma_{ii'} k_i k_{i'}} X_1^{k_1} X_2^{k_2}
\dots X_n^{k_n}
\end{align*}
if $q= \mathrm e^{\pi \mathrm i \hbar}$. In particular, this makes
the invariance under reindexing immediate.
\end{rem}

\begin{prop}
\label{prop:QandHinvariant} For any two ideal triangulations
$\lambda$, $\lambda' \in \Lambda(S)$, the coordinate change
isomorphism
    $\Phi_{\lambda\lambda'}^q:
\widehat{\mathcal{T}}^q_{\lambda'} \rightarrow
\widehat{\mathcal{T}}^q_{\lambda}$ sends $P_j \in
\widehat{\mathcal{T}}^q_{\lambda'}$ to $P_j \in
\widehat{\mathcal{T}}^q_{\lambda}$, and $H \in
\widehat{\mathcal{T}}^q_{\lambda'}$ to $H \in
\widehat{\mathcal{T}}^q_{\lambda}$.
\end{prop}

\begin{proof} Because of
Proposition~\ref{prop:InvarReindex}, it suffices to check this
when $\lambda'$ and $\lambda$ differ only by a diagonal exchange.

We will verify the property case-by-case, according to the type of
the diagonal exchange.  To help distinguishing the quantities
associated to $\lambda$  from those associated to $\lambda'$, we
will label with primes any data associated to $\lambda'$. Thus,
$\mathcal T_{\lambda'}^q$ is generated by the $X_i'$, and $P_j'
\in \mathcal T_{\lambda'}^q$ is the element associated to the
$j$--th puncture.

The diagonal exchange specifies two distinct triangle components
of $S-\lambda$. Let $\overline  Q$ be the closure of these two
triangles in $S$. Abstractly, $\overline  Q$ is obtained from the
square $Q$ where the diagonal exchange takes place by identifying
some of its sides according to the case we are considering. Note
that distinct vertices of $\overline Q$ can correspond to the same
puncture $v_j$ of $S$.

To simplify the notation (by saving the letter $j$), we focus
attention on the element $P_1$ associated to the first puncture
$v_1$. Since (the argument of) Proposition~\ref{prop:InvarReindex}
guarantees that the order in which we write its generators does
not matter, we can write $ P_1 = q^{-\sum_{k<l} \sigma_{i_k
i_l}^\lambda}X_{i_1} X_{i_2}\dots X_{i_m} $ where $X_{i_1}$,
$X_{i_2}$, \dots $X_{i_m}$ occur in this order as one goes
counterclockwise around the puncture $v_1$. Similarly, $ P_1' =
q^{-\sum_{k<l} \sigma_{i_k i_l}^{\lambda'}} X_{i_1}' X_{i_2}'\dots
X_{i_m}' $.

After a possible cyclic permutation of the $X_{i_u}$, those
touching the vertex $v$ of $\overline Q$ correspond to an interval
$a\leq k \leq b$. Let the \emph{contribution} of $v$ to $P_1$ be
$$
c_v(P_1) = q^{-\sum_{a\leq k<l\leq b} \sigma_{i_k
i_l}^\lambda}X_{i_a} X_{i_{a+1}}\dots X_{i_b}.
$$
The contribution $c_v(P_1')$ of $v$ to $P_1'$ is similarly
defined.

The element $P_1$ may contain the contribution of several vertices
of $\overline Q$. It is immediate that the part of $P_1$ which
does not come from these contributions is invariant under the
diagonal exchange. Therefore, we only need to show that
$\Phi_{\lambda\lambda'}^q$ sends $c_v(P_1') $ to $c_v(P_1) $ for
every vertex $v$ of $\overline Q$.

    We
label the data of the diagonal exchange according to
Figure~\ref{fig:DiagExch} and to the cases of
Proposition~\ref{prop:QuantumDiagExch}.

\medskip
\noindent\textbf{Case 1:}  the edges $\lambda_j$, $\lambda_k$,
$\lambda_l$, $\lambda_m$  are distinct (embedded diagonal
exchange).

Because of the symmetries of the figure and of the fact that this
case is its own inverse, it suffices to show this for an arbitrary
corner of $\overline Q=Q$, for instance the upper right corner. If
$v$ contributes to $P_1$ (and $P_1'$),
$$ c_v(P_1')=q^{-1}X_j'X_k' \quad\textrm{and}\quad
c_v( P_1) = q^{-2} X_jX_iX_k.
$$
It follows that
\begin{align*}
\Phi_{\lambda\lambda'}^q (c_v(P_1') )
&=q^{-1}(1+qX_i)X_j(1+qX_i^{-1})^{-1}X_k \\{} & =  q^{-1}X_j\,
(1+q^{-1}X_i)(1+qX_i^{-1})^{-1}X_k
\\
&= q^{-1} X_j\,q^{-1}X_iX_k= q^{-2} X_jX_iX_k= c_v(P_1).
\end{align*}

\medskip
\noindent\textbf{Cases 2 (and Case 3):}  $\lambda_j$ is identified
with $\lambda_k$, and $\lambda_l$ is distinct from $\lambda_m$.

We restrict to Case~2 since Case 3 is the inverse of Case 2.

In these cases, $\overline Q$ has three distinct vertices. The
first one corresponds to the upper left and lower right vertices
of $Q$. If this vertex $v$ contributes to $P_1$,
$$ c_v(P_1')=q^{-1} X_m'X_i' X_j' X_i' X_l'
\quad\textrm{and}\quad c_v( P_1) =  X_mX_jX_l.
$$
Then,
\begin{align*}
\Phi_{\lambda\lambda'}^q (c_v(P_1') ) &=q^{-1}
(1+qX_i^{-1})^{-1}X_m X_i^{-1}X_i X_j X_i^{-1} (1+qX_i) X_l
\\
&= q^{-1} X_m X_j(1+q^{-1}X_i^{-1})^{-1} X_i^{-1} (1+qX_i) X_l
\\
&= X_m X_j X_l = c_v( P_1)
\end{align*}

The second vertex of $\overline Q$ corresponds to the upper right
vertex of $Q$. For this vertex $v$,
$$ c_v(P_1')= X_j'\quad
\textrm{and}\quad c_v( P_1) =  X_i X_j.
$$
Then
$$
\Phi_{\lambda\lambda'}^q (c_v(P_1') ) = X_i X_j = c_v( P_1).
$$

The argument for the third vertex of $\overline Q$, corresponding
to the lower left vertex of $Q$, is identical to that of Case~1.

\medskip
\noindent\textbf{Cases 4 (and Case 5):}  $\lambda_j$ is identified
with $\lambda_l$, and $\lambda_k$ is distinct from $\lambda_m$.

We again can restrict attention to Case~4 since Case~5 is the
inverse of Case~4.

There are two vertices in $\overline Q$. This first one
corresponds to the two left corners of $Q$. If this vertex
contributes to $P_1$,
$$ c_v(P_1')= q^2 X_m' X_i' X_j' X_m\quad
\textrm{and}\quad c_v( P_1) =  q^2 X_m X_j X_i X_m.
$$
Then,
\begin{align*}
\Phi_{\lambda\lambda'}^q (c_v(P_1') ) &=q^2 (1+qX_i^{-1})^{-1}X_m
X_i^{-1} (1 +qX_i) (1+q^3 X_i) X_j (1+qX_i^{-1})^{-1}X_m
\\
&=q^2 X_m (1+q^{-1}X_i^{-1})^{-1} X_i^{-1} (1 +qX_i) (1+q^3 X_i)
(1+q^{-3}X_i^{-1})^{-1}X_jX_m
\\
&= q^2 X_m q^4 X_i X_j X_m = q^2 X_m X_j X_i X_m = c_v( P_1).
\end{align*}
The argument is similar for the second vertex of $\overline Q$,
corresponding to the right corners of $Q$.

\medskip
\noindent\textbf{Cases 6-8:} The arguments is similar to the above
cases.

\medskip
The same arguments apply to show that
$\Phi_{\lambda\lambda'}^q(H') = H$.
\end{proof}

Proposition~\ref{prop:QandHinvariant} shows that $H$ and the $P_j$
are well-defined elements of the quantum enhanced Teichm\"uller
space $\widetilde{\mathcal{T}}^q_S$. It is not too hard to see
that these elements are central. It is proved in \cite{BonLiu}
that they generate the center of $\widetilde{\mathcal{T}}^q_S$.

By analogy with the non-quantum case, we define the \emph{quantum
cusped Teichm\"uller space} $\mathcal{CT}^q_S$ as the quotient
algebra
$$
\mathcal{CT}^q_S= \widetilde{\mathcal{T}}^q_S/I
$$
where $I$ is the 2--sided ideal generated by $H-\mathbf{1}$ and
the $P_j-\mathbf{1}$. As indicated at the beginning of this
section, the algebra $\mathrm{Rat}\,\mathcal {CT}(S)$ on the
cusped Teichm\"uller space $\mathcal {CT}(S)$ is naturally
isomorphic to $\mathcal{CT}^1_S$ when $q=1$.

In \cite{BonLiu}, we define invariants of a diffeomorphism
$\varphi$ of $S$ by considering certain representations $\rho$ of
the polynomial core (defined in that paper) of the quantum
enhanced  Teichm\"uller space $\widetilde{\mathcal T}_S^q$. These
representations are associated to the (unique) hyperbolic metric
on the mapping torus $M_\varphi = S\times \mathbb R/\sim$, where
$\sim$ identifies $(x,t)$ to $ (\varphi^n(x), t +n)$. In
particular, a property of $\rho$ is that $\rho(H)=\mathrm{Id}$ and
$\rho(P_j)=\mathrm{Id}$. In other words, this representation
$\rho$ is actually a representation of the polynomial core of the
quantum cusped Teichm\"uller space.

\end{document}